\input amstex.tex
\documentstyle{amsppt}
\magnification=\magstep1
\vsize 8.5truein
\hsize 6truein

\define\flow{\left(\bold{M},\{S^t\}_{t\in\Bbb R},\mu\right)}

\noindent
March 8, 2007

\bigskip \bigskip

\heading
An Effective Contraction Estimate in the Stable Subspaces of Phase Points in 
Hard Ball Systems
\endheading

\bigskip \bigskip

\centerline{{\bf N\'andor Sim\'anyi}
\footnote{Research supported by the National Science Foundation, 
grant DMS-0457168.}}

\bigskip \bigskip

\centerline{University of Alabama at Birmingham}
\centerline{Department of Mathematics}
\centerline{Campbell Hall, Birmingham, AL 35294 U.S.A.}
\centerline{E-mail: simanyi\@math.uab.edu}

\bigskip \bigskip

\hbox{\centerline{\vbox{\hsize 8cm{\bf Abstract.} In this paper we prove the
      following result, useful and often needed in the study of the ergodic
      properties of hard ball systems: In any hard ball system, for any phase
      point $x_0\in\bold M\setminus\partial\bold M$ with a non-singular
      forward orbit $S^{(0,\infty)}x_0$ and with infinitely many consecutive,
      connected collision graphs on $S^{(0,\infty)}x_0$, and for any number
      $L>0$ one can always find a time $t>0$ and a non-zero tangent vector
      $(\delta q_0,\delta v_0)\in E^s(x_{0})$ with $\frac{||(\delta
        q_t,\,\delta v_t)||}{||(\delta q_0,\,\delta v_0)||}<L^{-1},$ where
      $(\delta q_t,\delta v_t)=DS^t(\delta q_0,\delta v_0)\in E^s(x_t)$,
      $x_t=S^tx_0$. Of course, the Multiplicative Ergodic Theorem of Oseledets
      provides a much stronger conclusion, but at the expense of an
      unspecified zero-measured exceptional set of phase points, and this is
      not sufficient in the sophisticated studies the ergodic properties of
      such flows. Here the exceptional set of phase points is a dynamically
      characterized set, so that it suffices for the proofs showing how global
      ergodicity follows from the local one.}}}

\bigskip

\noindent
Primary subject classification: 37D50

\medskip

\noindent
Secondary subject classification: 34D05

\bigskip \bigskip

\heading
\S1. Introduction
\endheading

\bigskip

Ever so often it happens that, while studying the ergodic properties of
semi-dispersive billiard flows, one needs the following result:

\medskip

\subheading{Theorem} For any phase point 
$x_0\in\bold M\setminus\partial\bold M$ with a non-singular
forward orbit $S^{(0,\infty)}x_0$ and with infinitely many consecutive,
connected collision graphs on $S^{(0,\infty)}x_0$, and for any
number $L>0$ one can always find a time $t>0$ and a non-zero tangent vector 
$(\delta q_0,\delta v_0)\in E^s(x_{0})$ with

$$
\frac{||(\delta q_t,\,\delta v_t)||}{||(\delta q_0,\,\delta v_0)||}<L^{-1},
$$
where $(\delta q_t,\delta v_t)=DS^t(\delta q_0,\delta v_0)\in E^s(x_t)$,
$x_t=S^tx_0$.

\medskip

The last time when the above result was badly needed, at least in my own
research, was a certain step in the conditional proof of the Boltzmann-Sinai
Ergodic Hypothesis in [Sim(2006)].

\medskip

We note that, according to Oseledets' Multiplicative Ergodic Theorem and the
full hyperbolicity of every hard ball system [Sim(2002)], for almost every
phase point $x_0$ (with respect to the natural invariant measure, that is, the
Liouville measure) every vector $w\in E^s(x_0)$ contracts exponentially fast
as $t\to\infty$. The point in the present result is that it holds true for
even more generic points $x_0$: for all points apart from a so called slim
subset of the phase space, see Theorem 5.1 in [Sim(1992)-I]. What is more,
even if we have a smooth, codimension-one, exceptional manifold $J$ (capable
of separating two distinct, open ergodic components of the flow) possessing a
normal vector field $n(x)=(z(x),\,w(x))$ with $Q(n(x))=\langle
z(x),\,w(x)\rangle<0$ ($x\in J$, see the beginning of \S3 in [Sim(2006)]),
almost every phase point $x_0$ of the hypersurface $J$ will satisfy the
hypotheses of the above theorem. This is indeed so, since the proof of Theorem
6.1 of [Sim(1992)-I] works without any essential change not only for singular
phase points, but also for the points of the considered exceptional manifold
$J$. The only important ingredient of that proof is the transversality of the
spaces $E^s(x)$ to $J$ for all $x\in J$. According to that result, typical
phase points $x\in J$ (with respect to the hypersurface measure of $J$) indeed
enjoy the above property of having infinitely many consecutive, connected
collision graphs on their forward orbit $S^{(0,\infty)}x_0$.

\medskip

We also note that, by slightly perturbing the shrinking vector $(\delta
q_0,\delta v_0)\in E^s(x_{0})$ {\it inside the stable space} $E^s(x_{0})$, we
can achieve that this vector be transversal to the manifold $J$. This is
another point where we use the fact that the stable spaces $E^s(x)$ ($x\in J$)
are transversal to $J$, an immediate corollary of the negativity assumption
on the infinitesimal Lyapunov form $Q(n(x))$.

\medskip

It is worth noting here that the main phenomenon that makes the assertion of
the theorem substantial (and the proof non-trivial) is the possibility of very
long free flights (i. e. orbit segments without collision), or at least the
possibility of long segments on $S^{(0,\infty)}x_0$ in which all occuring
collisions have very small relative velocities for the colliding balls. This
phenomenon ``flattens out'' the stable manifold $\gamma^s(x_t)$ of
$x_t=S^tx_0$ and, as a result, has the potential to make the contraction
coefficient only slightly smaller than $1$. The main part of the proof of the
theorem is to show that the above phenomenon actually does not occur, at least
when $S^{(0,\infty)}x_0$ has infinitely many consecutive, connected
collision graphs (Corollary 3.12). It is just the proof of this corollary,
more precisely, the proof of Proposition 3.9, that utilizes the assumption on
the infinitely many connected collision graphs.

\medskip

Finally, we note that the unstable version of the above theorem, claiming the
arbitrarily big expansions (as $t\to \infty$) of the forward images of vectors
$w\in E^u(x_0)$ is obviously true: It is easy to see that any tangent vector
$w=(\delta q,\,\delta v)\in\Cal T_{x_0}\bold M$ with $\langle\delta q,\,\delta
v\rangle>0$ expands at least linearly in time as $t\to\infty$, even without
collisions, see Proposition 3.5 below.

\bigskip \bigskip

\heading
\S2. Prerequisites
\endheading

\bigskip

Consider the $\nu$-dimensional ($\nu\ge2$), standard, flat torus
$\Bbb T^\nu=\Bbb R^\nu/\Bbb Z^\nu$ as the vessel containing 
$N$ ($\ge2$) hard balls (spheres) $B_1,\dots,B_N$ with positive masses 
$m_1,\dots,m_N$ and (just for simplicity) common radius $r>0$. We always
assume that the radius $r>0$ is not too big, so
that even the interior of the arising configuration space $\bold Q$ (or, 
equivalently, the phase space) is connected. Denote the center of the ball
$B_i$ by $q_i\in\Bbb T^\nu$, and let $v_i=\dot q_i$ be the velocity of the
$i$-th particle. We investigate the uniform motion of the balls
$B_1,\dots,B_N$ inside the container $\Bbb T^\nu$ with half a unit of total 
kinetic energy: $E=\dfrac{1}{2}\sum_{i=1}^N m_i||v_i||^2=\dfrac{1}{2}$.
We assume that the collisions between balls are perfectly elastic. Since
--- beside the kinetic energy $E$ --- the total momentum
$I=\sum_{i=1}^N m_iv_i\in\Bbb R^\nu$ is also a trivial first integral of the
motion, we make the standard reduction $I=0$. Due to the apparent translation
invariance of the arising dynamical system, we factorize the configuration
space with respect to uniform spatial translations as follows:
$(q_1,\dots,q_N)\sim(q_1+a,\dots,q_N+a)$ for all translation vectors
$a\in\Bbb T^\nu$. The configuration space $\bold Q$ of the arising flow
is then the factor torus
$\left(\left(\Bbb T^\nu\right)^N/\sim\right)\cong\Bbb T^{\nu(N-1)}$
minus the cylinders

$$
C_{i,j}=\left\{(q_1,\dots,q_N)\in\Bbb T^{\nu(N-1)}\colon\;
\text{dist}(q_i,q_j)<2r \right\}
$$
($1\le i<j\le N$) corresponding to the forbidden overlap between the $i$-th
and $j$-th spheres. Then it is easy to see that the compound 
configuration point

$$
q=(q_1,\dots,q_N)\in\bold Q=\Bbb T^{\nu(N-1)}\setminus
\bigcup_{1\le i<j\le N}C_{i,j}
$$
moves in $\bold Q$ uniformly with unit speed and bounces back from the
boundaries $\partial C_{i,j}$ of the cylinders $C_{i,j}$ according to the
classical law of geometric optics: the angle of reflection equals the angle of
incidence. More precisely: the post-collision velocity $v^+$ can be obtained
from the pre-collision velocity $v^-$ by the orthogonal reflection across the
tangent hyperplane of the boundary $\partial\bold Q$ at the point of collision.
Here we must emphasize that the phrase ``orthogonal'' should be understood 
with respect to the natural Riemannian metric (the kinetic energy)
$||dq||^2=\sum_{i=1}^N m_i||dq_i||^2$ in the configuration space $\bold Q$.
For the normalized Liouville measure $\mu$ of the arising flow
$\{S^t\}$ we obviously have $d\mu=\text{const}\cdot dq\cdot dv$, where
$dq$ is the Riemannian volume in $\bold Q$ induced by the above metric, 
and $dv$ is the surface measure (determined by the restriction of the
Riemannian metric above) on the unit sphere of compound velocities

$$
\Bbb S^{\nu(N-1)-1}=\left\{(v_1,\dots,v_N)\in\left(\Bbb R^\nu\right)^N\colon\;
\sum_{i=1}^N m_iv_i=0 \text{ and } \sum_{i=1}^N m_i||v_i||^2=1 \right\}.
$$
The phase space $\bold M$ of the flow $\{S^t\}$ is the unit tangent bundle
$\bold Q\times\Bbb S^{d-1}$ of the configuration space $\bold Q$. (We will 
always use the shorthand notation $d=\nu(N-1)$ for the dimension of the 
billiard table $\bold Q$.) We must, however, note here that at the boundary
$\partial\bold Q$ of $\bold Q$ one has to glue together the pre-collision and
post-collision velocities in order to form the phase space $\bold M$, so
$\bold M$ is equal to the unit tangent bundle $\bold Q\times\Bbb S^{d-1}$
modulo this identification.

A bit more detailed definition of hard ball systems with arbitrary masses,
as well as their role in the family of cylindric billiards, can be found in
\S4 of [S-Sz(2000)] and in \S1 of [S-Sz(1999)]. We denote the
arising flow by $\flow$.

\medskip

\subheading{Collision graphs} Let $S^{[a,b]}x$ be a nonsingular, finite
trajectory segment with the collisions $\sigma_1,\dots,\sigma_n$
listed in time order. 
(Each $\sigma_k$ is an unordered pair $(i,j)$ of different labels
$i,j\in\{1,2,\dots,N\}$.) The graph $\Cal G=(\Cal V,\Cal E)$ with vertex set
$\Cal V=\{1,2,\dots,N\}$ and set of edges $\Cal E=\{\sigma_1,\dots,\sigma_n\}$ 
is called the {\it collision graph}
of the orbit segment $S^{[a,b]}x$. For a given positive number $C$, the
collision graph $\Cal G=(\Cal V,\Cal E)$ of the orbit segment $S^{[a,b]}x$
will be called {\it $C$-rich} if $\Cal G$ contains at least $C$ connected,
consecutive (i. e. following one after the other in time, according to the
time-ordering given by the trajectory segment $S^{[a,b]}x$) subgraphs. 

\medskip

\subheading{No accumulation (of collisions) in finite time} 
By the results of Vaserstein [V(1979)], Galperin [G(1981)] and
Burago-Ferleger-Kononenko [B-F-K(1998)], in any semi-dis\-per\-sive 
billiard flow there can only be finitely many 
collisions in finite time intervals, see Theorem 1 in [B-F-K(1998)]. 
Thus, the dynamics is well defined as long as the trajectory does not hit 
more than one boundary components at the same time.

\medskip

Finally, for any phase point $x\in\bold M\setminus\partial\bold M$ with a
non-singular forward orbit $S^{(0,\infty)}x$ (and with at least one collision,
hence infinitely many collisions on it) we define the stable subspace 
$E^s(x)\subset\Cal T_x\bold M$ of $x$ as

$$
E^s(x)=\left\{(\delta q,\,\delta v)\in\Cal T_x\bold M\big|\; 
\delta v=-B(x)[\delta q],\,\, \langle\delta q,\,\delta v\rangle<0\right\}
\cup\left\{(0,\,0)\right\},
$$
where the symmetric, non-negative operator $B(x)$ (acting on the tangent space
of $\bold Q$ at the footpoint $q$, where $x=(q,v)$) is defined by the
continued fraction expansion introduced by Sinai in [Sin(1979)], see also
[Ch(1982)] or (2.4) in [K-S-Sz(1990)-I]. It is a well known fact that
$E^s(x)$ is the tangent space of the local stable manifold $\gamma^s(x)$, if
the latter object exists.

\medskip

For any phase point $x\in\bold M\setminus\partial\bold M$ with a non-singular
backward orbit $S^{(-\infty,0)}x$ (and with at least one collision on it) the
unstable space $E^u(x)$ of $x$ is defined as $E^s(-x)$, where 
$-x=(q,\,-v)$ for $x=(q,\,v)$.

\bigskip\bigskip

\heading
\S3. Expansion and Contraction Rate Estimates \\
Proof of Theorem
\endheading

\bigskip\bigskip

We would like to get a useful lower estimate for the expansion of a tangent
vector $(\delta q_0,\delta v_0)\in\Cal T_{x_0}\bold M$ with positive
infinitesimal Lyapunov function $Q(\delta q_0,\delta v_0)=\langle\delta
q_0,\delta v_0\rangle$. The expression $\langle\delta q_0,\delta v_0\rangle$
is the scalar product in $\Bbb R^d$ defined via the mass (or kinetic energy)
metric, see \S2. It is also called the infinitesimal Lyapunov function
associated with the tangent vector $(\delta q_0,\delta v_0)$, see [K-B(1994)],
or part A.4 of the Appendix in [Ch(1994)], or \S7 of [Sim(2003)]. For a
detailed exposition of the relationship between the quadratic form $Q(\,.\,)$,
the relevant symplectic geometry of the Hamiltonian system and the dynamics,
please also see [L-W(1995)].

\medskip

\subheading{Note} The original idea of using infinitesimal Lyapunov exponents
to measure the expansion rate of codimension-one submanifolds in the phase
space of semi-dispersive billiards came from N. Chernov back in the late '80s.
These ideas have been explored in detail and further developed by him and
myself in recent personal communications, so that we obtained at least linear
(but uniform!)  expansion rates for such submanifolds with negative
infinitesimal Lyapunov forms for their normal vector. These results are
presented in our recent joint paper [Ch-Sim(2006)]. Also, closely related to
the above said, the following ideas (to estimate the expansion rates of tangent
vectors from below) are derived from the thoughts being published in 
[Ch-Sim(2006)].

\medskip

Denote by $(\delta q_t,\delta v_t)=(DS^t)(\delta q_0,\delta v_0)$
the image of the tangent vector $(\delta q_0,\delta v_0)$ under the
linearization $DS^t$ of the map $S^t$, $t\ge0$. (We assume that the base phase
point $x_0$ --- for which $(\delta q_0,\delta v_0)\in\Cal T_{x_0}\bold M$ ---
has a non-singular forward orbit.) The time-evolution 
$(\delta q_{t_1},\delta v_{t_1})\mapsto(\delta q_{t_2},\delta v_{t_2})$
($0\le t_1<t_2$) on a collision free segment $S^{[t_1,t_2]}x_0$ is described
by the equations

$$
\aligned
\delta v_{t_2}&=\delta v_{t_1}, \\
\delta q_{t_2}&=\delta q_{t_1}+(t_2-t_1)\delta v_{t_1}.
\endaligned
\tag 3.1
$$
Correspondingly, the change 
$Q(\delta q_{t_1},\delta v_{t_1})\mapsto Q(\delta q_{t_2},\delta v_{t_2})$
in the infinitesimal Lyapunov function $Q(\,.\,)$ on the collision free orbit
segment $S^{[t_1,t_2]}x_0$ is

$$
Q(\delta q_{t_2},\delta v_{t_2})=Q(\delta q_{t_1},\delta v_{t_1})
+(t_2-t_1)||\delta v_{t_1}||^2,
\tag 3.2
$$
thus $Q(\,.\,)$ steadily increases between collisions. 

The passage
$(\delta q_t^-,\delta v_t^-)\mapsto(\delta q_t^+,\delta v_t^+)$ through a
reflection (i. e. when $x_t=S^tx_0\in\partial\bold M$) is given by Lemma 2 of
[Sin(1979)] or formula (2) in \S3 of [S-Ch(1987)]:

$$
\aligned
\delta q_t^+&=R\delta q_t^-, \\
\delta v_t^+&=R\delta v_t^-+2\cos\phi RV^*KV\delta q_t^-,
\endaligned
\tag 3.3
$$
where the operator $R:\, \Cal T\bold Q\to\Cal T\bold Q$ is the orthogonal 
reflection (with respect to the mass metric) across the tangent hyperplane
$\Cal T_{q_t}\partial\bold Q$ of the boundary $\partial\bold Q$ at the
configuration component $q_t$ of $x_t=(q_t,v_t^\pm)$, 
$V:\, (v_t^-)^\perp\to\Cal T_{q_t}\partial\bold Q$ is the $v_t^-$-parallel
projection of the orthocomplement hyperplane $(v_t^-)^\perp$ onto
$\Cal T_{q_t}\partial\bold Q$, 
$V^*:\,\Cal T_{q_t}\partial\bold Q\to(v_t^-)^\perp$ is the adjoint of $V$
(i. e. the $\nu(q_t)$-parallel projection of $\Cal T_{q_t}\partial\bold Q$
onto $(v_t^-)^\perp$, where $\nu(q_t)$ is the inner normal vector of 
$\partial\bold Q$ at $q_t\in\partial\bold Q$), 
$K:\, \Cal T_{q_t}\partial\bold Q\to\Cal T_{q_t}\partial\bold Q$ is the second
fundamental form of the boundary $\partial\bold Q$ at $q_t$ (with respect to 
the field $\nu(q)$ of inner unit normal vectors of $\partial\bold Q$) and,
finally, $\cos\phi=\langle\nu(q_t),\, v_t^+\rangle>0$ is the cosine of the 
angle $\phi$ ($0\le\phi<\pi/2$) subtended by $v_t^+$ and $\nu(q_t)$. Regarding
formulas (3.3), please see the last displayed formula in \S1 of 
[S-Ch(1987)] or (i)--(ii) in Proposition 2.3 of [K-S-Sz(1990)-I]. The
instanteneous change in the infinitesimal Lyapunov function
$Q(\delta q_t,\delta v_t)$ caused by the reflection at time $t>0$ is easily
derived from (3.3):

$$
\aligned
Q(\delta q_t^+,\delta v_t^+)&=Q(\delta q_t^-,\delta v_t^-)+2\cos\phi\langle
V\delta q_t^-,\, KV\delta q_t^-\rangle \\
&\ge Q(\delta q_t^-,\delta v_t^-).
\endaligned
\tag 3.4
$$ 
In the last inequality we used the fact that the operator $K$ is positive
semi-definite, i. e. the billiard is semi-dispersive.

We are primarily interested in getting useful lower estimates for the expansion
rate $||\delta q_t||/||\delta q_0||$. The needed result is

\medskip

\subheading{Proposition 3.5} Use all the notations above, and assume that

$$
\langle\delta q_0,\, \delta v_0\rangle/||\delta q_0||^2\ge c_0>0.
$$
We claim that $||\delta q_t||/||\delta q_0||\ge 1+c_0t$ for all $t\ge0$.

\medskip

\subheading{Proof} Clearly, the function $||\delta q_t||$ of $t$ is continuous
for all $t\ge0$ and continuously differentiable between collisions. According 
to (3.1), $\frac{d}{dt}\delta q_t=\delta v_t$, so

$$
\frac{d}{dt}||\delta q_t||^2=2\langle\delta q_t, \delta v_t\rangle.
\tag 3.6
$$

Observe that not only the positive valued function 
$Q(\delta q_t,\delta v_t)=\langle\delta q_t,\delta v_t\rangle$ is 
nondecreasing in $t$ by (3.2) and (3.4), but the quantity 
$\langle\delta q_t,\delta v_t\rangle/||\delta q_t||$ is nondecreasing in $t$,
as well. The reason is that
$\langle\delta q_t,\delta v_t\rangle/||\delta q_t||=||\delta v_t||\cos\alpha_t$
($\alpha_t$ being the acute angle subtended by $\delta q_t$ and $\delta v_t$),
and between collisions the quantity $||\delta v_t||$ is unchanged, while the
acute angle $\alpha_t$ decreases, according to the time-evolution equations
(3.1). Finally, we should keep in mind that at a collision the norm
$||\delta q_t||$ does not change, while $\langle\delta q_t,\delta v_t\rangle$
cannot decrease, see (3.4). Thus we obtain the inequalities

$$
\langle\delta q_t,\delta v_t\rangle/||\delta q_t||\ge
\langle\delta q_0,\delta v_0\rangle/||\delta q_0||\ge c_0||\delta q_0||,
$$
so

$$
\frac{d}{dt}||\delta q_t||^2=2||\delta q_t||\frac{d}{dt}||\delta q_t||
=2\langle\delta q_t,\delta v_t\rangle
\ge 2c_0||\delta q_0||\cdot||\delta q_t||
$$
by (3.6). This means that 
$\frac{d}{dt}||\delta q_t||\ge c_0||\delta q_0||$, so
$||\delta q_t||\ge||\delta q_0||(1+c_0t)$, proving the proposition. \qed

\medskip

Next we need an effective lower estimation $c_0$ for the curvature
$\langle\delta q_0,\, \delta v_0\rangle/||\delta q_0||^2$ of the trajectory
bundle:

\medskip

\subheading{Lemma 3.7} Assume that the perturbation 
$(\delta q_0^-,\,\delta v_0^-)\in\Cal T_{x_0}\bold M$ (as in Proposition 3.5)
is being performed at time zero right before a collision, say,
$\sigma_0=(1,\,2)$ taking place at that time. Select the tangent vector
$(\delta q_0^-,\,\delta v_0^-)$ in such a specific way that $\delta
q_0^-=(m_2w,-m_1w,0,0,\dots,0)$ with a nonzero vector $w\in\Bbb R^\nu$,
$\langle w,v_1^--v_2^-\rangle=0$. This scalar product equation is exactly the
condition that guarantees that $\delta q_0^-$ be orthogonal to the velocity
component $v^-=(v_1^-,v_2^-,\dots,v_N^-)$ of $x_0=(q,v^-)$. The next, though
crucial requirement is that $w$ should be selected from the two-dimensional
plane spanned by $v_1^--v_2^-$ and $q_1-q_2$ (with $||q_1-q_2||=2r$) in $\Bbb
R^\nu$. The purpose of this condition is to avoid the unwanted phenomenon of
``astigmatism'' in our billiard system, discovered first by Bunimovich and
Rehacek in [B-R(1997)] and [B-R(1998)]. Later on the phenomenon of astigmatism
gathered further prominence in the paper [B-Ch-Sz-T(2002)] as the main driving
mechanism behind the wild non-differentiability of the singularity manifolds
(at their boundaries) in hard ball systems in dimensions bigger than
$2$. Finally, the last requirement is that the velocity component $\delta
v_0^-$ (right before the collision $(1,2)$) is chosen in such a way that the
tangent vector $(\delta q_0^-,\,\delta v_0^-)$ belongs to the unstable space
$E^u(x_0)$ of $x_0$. This can be done, indeed, by taking $\delta
v_0^-=B^u(x_0)[\delta q_0^-]$, where $B^u(x_0)$ the curvature operator of the
unstable manifold of $x_0$ at $x_0$, right before the collision $(1,2)$ taking
place at time zero. Note that the argument $\delta q_0^-$ of $B^u(x_0)$
clearly belongs to the positive subspace of the non-negative operator
$B^u(x_0)$.

We claim that

$$
\frac{\langle\delta q_0^+,\delta v_0^+\rangle}{||\delta q_0||^2}\ge
\frac{||v_1-v_2||}{r\cos\phi_0}\ge\frac{||v_1-v_2||}{r}
\tag 3.8
$$
for the post-collision tangent vector $(\delta q_0^+,\delta v_0^+)$, where
$\phi_0$ is the acute angle subtended by $v_1^+-v_2^+$ and the outer normal
vector of the sphere $\left\{y\in\Bbb R^\nu\big|\; ||y||=2r\right\}$ at the 
point $y=q_1-q_2$. Note that in (3.8) there is no need to use $+$ or $-$ in 
$||\delta q_0||^2$ or $||v_1-v_2||$, for $||\delta q_0^-||=||\delta q_0^+||$,
$||v_1^--v_2^-||=||v_1^+-v_2^+||$.

\medskip

\subheading{Proof} The proof of the equation in (3.8) is a simple, elementary
geometric argument in the plane spanned by $v_1^--v_2^-$ and $q_1-q_2$, so we
omit it. We only note that the outgoing relative velocity $v_1^+-v_2^+$ is
obtained from the pre-collision relative velocity $v_1^--v_2^-$ by reflecting
the latter one across the tangent hyperplane of the sphere 
$\left\{y\in\Bbb R^\nu\big|\; ||y||=2r\right\}$ at the point $y=q_1-q_2$.
It is a useful advice, though, to prove the first inequality of (3.8) in the
case $\delta v_0^-=0$ first (this is the elementary geometry exercise), then
observe that this inequality can only be further improved when we replace 
$\delta v_0^-=0$ with $\delta v_0^-=B^u(x_0)[\delta q_0^-]$. \qed

\medskip

The previous lemma shows that, in order to get useful lower estimates for the
``curvature'' $\langle\delta q,\delta v\rangle/||\delta q||^2$ of the
trajectory bundle, it is necessary
(and sufficient) to find collisions $\sigma=(i,j)$ on the orbit of a given 
point $x_0\in\bold M$ with a ``relatively big'' value of $||v_i-v_j||$.
Finding such collisions will be based upon the following result:

\medskip

\subheading{Proposition 3.9} Consider orbit segments $S^{[0,T]}x_0$ of 
$N$-ball systems with masses $m_1,m_2,\dots,m_N$ in $\Bbb T^\nu$ (or in
$\Bbb R^\nu$) with collision sequences 
$\Sigma=(\sigma_1,\sigma_2,\dots,\sigma_n)$ corresponding to connected 
collision graphs. (Now the kinetic energy is not necessarily normalized, and 
the total momentum $\sum_{i=1}^N m_iv_i$ may be different from zero.) We claim
that there exists a positive-valued function $f(a;m_1,m_2,\dots,m_N)$ 
($a>0$, $f$ is independent of the orbit segments $S^{[0,T]}x_0$) with the
following two properties:

\medskip

(1) If $||v_i(t_l)-v_j(t_l)||\le a$ for all collisions $\sigma_l=(i,j)$
($1\le l\le n$, $t_l$ is the time of $\sigma_l$) of some trajectory segment 
$S^{[0,T]}x_0$ with a symbolic collision sequence 
$\Sigma=(\sigma_1,\sigma_2,\dots,\sigma_n)$ corresponding to a connected 
collision graph, then the norm $||v_{i'}(t)-v_{j'}(t)||$ of any relative
velocity at any time $t\in\Bbb R$ is at most $f(a;m_1,\dots,m_N)$;

(2) $\lim_{a\to 0} f(a;m_1,\dots,m_N)=0$ for any $(m_1,\dots,m_N)$.

\medskip

\subheading{Proof} We begin with

\medskip

\subheading{Lemma 3.10} Consider an $N$-ball system with masses $m_1,\dots,m_N$
(an $(m_1,\dots\allowmathbreak,m_N)$-system, for short) in $\Bbb T^\nu$ 
(or in $\Bbb R^\nu$).
Assume that the inequalities $||v_i(0)-v_j(0)||\le a$ hold true 
($1\le i<j\le N$) for all relative velocities at time zero. We claim that 

$$
||v_i(t)-v_j(t)||\le 2a\sqrt{\frac{M}{m}}
\tag 3.11
$$
for any pair $(i,j)$ and any time $t\in\Bbb R$, where $M=\sum_{i=1}^N m_i$ and 

$$
m=\min\left\{m_i|\; 1\le i\le N\right\}.
$$

\medskip

\subheading{Note} The estimate (3.11) is far from optimal, however, it
will be sufficient for our purposes.

\medskip

\subheading{Proof} The assumed inequalities directly imply that 
$||v'_i(0)||\le a$ ($1\le i\le N$) for the velocities $v'_i(0)$ measured at 
time zero in the baricentric reference system. Therefore, for the total kinetic
energy $E_0$ (measured in the baricentric system) we get the upper estimation
$E_0\le \frac{1}{2}Ma^2$, and this inequality remains true at any time $t$.
This means that all the inequalities $||v'_i(t)||^2\le\frac{M}{m_i}a^2$ hold
true for the baricentric velocities $v'_i(t)$ at any time $t$, so

$$
||v'_i(t)-v'_j(t)||\le a\sqrt{M}\left(m_i^{-1/2}+m_j^{-1/2}\right)\le
2a\sqrt{\frac{M}{m}},
$$
thus the inequalities

$$
||v_i(t)-v_j(t)||\le 2a\sqrt{\frac{M}{m}}
$$
hold true, as well. \qed

\medskip

\subheading{Proof of the proposition by induction on the number $N$}

\medskip

For $N=1$ we can take $f(a;m_1)=0$, and for $N=2$ the function
$f(a;m_1,m_2)=a$ is obviously a good choice for $f$. Let $N\ge 3$, and assume 
that the orbit segment $S^{[0,T]}x_0$ of an $(m_1,\dots,m_N)$-system fulfills
the conditions of the proposition. Let $\sigma_k=(i,j)$ be the collision in the
symbolic sequence $\Sigma_n=(\sigma_1,\dots,\sigma_n)$ of $S^{[0,T]}x_0$ with
the property that the collision graph of $\Sigma_k=(\sigma_1,\dots,\sigma_k)$
is connected, while the collision graph of 
$\Sigma_{k-1}=(\sigma_1,\dots,\sigma_{k-1})$ is still disconnected. Denote the
two connected components (as vertex sets) of $\Sigma_{k-1}$ by $C_1$ and
$C_2$, so that $i\in C_1$, $j\in C_2$, $C_1\cup C_2=\{1,2,\dots,N\}$, and
$C_1\cap C_2=\emptyset$. By the induction hypothesis and the condition of the
proposition, the norm of any relative velocity
$v_{i'}(t_k-0)-v_{j'}(t_k-0)$ right before the collision $\sigma_k$ (taking 
place at time $t_k$) is at most 
$a+f(a;\,\overline{C}_1)+f(a;\,\overline{C}_2)$, where $\overline{C}_l$ stands
for the collection of the masses of all particles in the component $C_l$,
$l=1,\,2$. Let $g(a;m_1,\dots,m_N)$ be the maximum of all possible sums

$$
a+f(a;\,\overline{D}_1)+f(a;\,\overline{D}_2),
$$
taken for all two-class partitions $(D_1,D_2)$ of the vertex set
$\{1,2,\dots,N\}$. According to the previous lemma, the function

$$
f(a;m_1,\dots,m_N):=2\sqrt{\frac{M}{m}}g(a;m_1,\dots,m_N)
$$
fulfills both requirements (1) and (2) of the proposition. \qed

\medskip

\subheading{Corollary 3.12} Consider the original $(m_1,\dots,m_N)$-system with
the standard normalizations $\sum_{i=1}^N m_iv_i=0$, 
$\frac{1}{2}\sum_{i=1}^N m_i||v_i||^2=\frac{1}{2}$. We claim that there exists
a threshold $G=G(m_1,\dots,m_N)>0$ (depending only on $N$, $m_1,\dots,m_N$)
with the following property:

In any orbit segment $S^{[0,T]}x_0$ of the $(m_1,\dots,m_N)$-system with the
standard normalizations and with a connected collision graph, one can always
find a collision $\sigma=(i,j)$, taking place at time $t$, so that
$||v_i(t)-v_j(t)||\ge G(m_1,\dots,m_N)$.

\medskip

\subheading{Proof} Indeed, we choose $G=G(m_1,\dots,m_N)>0$ so small that
$f(G;m_1,\dots,m_N)\allowmathbreak<M^{-1/2}$. 
Assume, contrary to 3.12, that the norm of any
relative velocity $v_i-v_j$ of any collision of $S^{[0,T]}x_0$ is less than the
above selected value of $G$. By the proposition, we have the inequalities
$||v_i(0)-v_j(0)||\le f(G;m_1,\dots,m_N)$ at time zero. The normalization
$\sum_{i=1}^N m_iv_i(0)=0$, with a simple convexity argument, implies that
$||v_i(0)||\le f(G;m_1,\dots,m_N)$ for all $i$, $1\le i\le N$, so the total
kinetic energy is at most
$\frac{1}{2}M\left[f(G;m_1,\dots,m_N)\right]^2<\frac{1}{2}$, a contradiction.
\qed

\medskip

\subheading{Corollary 3.13} For any phase point $x_0$ with a non-singular
backward trajectory $S^{(-\infty,0)}x_0$ and with infinitely many consecutive,
connected collision graphs on $S^{(-\infty,0)}x_0$, and for any number $L>0$
one can always find a time $-t<0$ and a non-zero tangent vector $(\delta
q_0,\delta v_0)\in E^u(x_{-t})$ ($x_{-t}=S^{-t}x_0$) with $||\delta
q_t||/||\delta q_0||>L$, where $(\delta q_t,\delta v_t)=DS^t(\delta q_0,\delta
v_0)\in E^u(x_0)$.

\medskip

\subheading{Proof} Indeed, select a number $t>0$ so big that
$1+\frac{t}{r}G(m_1,\dots,m_N)>L$ and $-t$ is the time of a collision (on the
orbit of $x_0$) with the relative velocity $v^-_i(-t)-v^-_j(-t)$, for which
$||v^-_i(-t)-v^-_j(-t)||\ge G(m_1,\dots,m_N)$. By Lemma 3.7 we can choose a
non-zero tangent vector $(\delta q^-_0,\,\delta v^-_0)\in E^u(x_{-t})$ right
before the collision at time $-t$ in such a way that the lower estimate

$$
\frac{\langle\delta q^+_0,\,\delta v^+_0\rangle}{||\delta q^+_0||^2}\ge
\frac{1}{r}G(m_1,\dots,m_N)
$$
holds true for the ``curvature''
$\langle\delta q^+_0,\,\delta v^+_0\rangle/||\delta q^+_0||^2$ associated
with the post-collision tangent vector $(\delta q^+_0,\,\delta v^+_0)$.
According to Proposition 3.5, we have the lower estimate

$$
\frac{||\delta q_t||}{||\delta q_0||}\ge 1+\frac{t}{r}G(m_1,\dots,m_N)>L
$$
for the $\delta q$-expansion rate between 
$(\delta q^-_0,\, \delta v^-_0)$ and
$(\delta q_t,\,\delta v_t)=DS^t(\delta q^-_0,\, \delta v^-_0)$. \qed

\medskip

We remind the reader that, according to the main result of [B-F-K(1998)],
there exists a number $\epsilon_0=\epsilon_0(m_1,\dots,m_N;\,r;\,\nu)>0$
and a large threshold $N_0=N_0(m_1,\dots,m_N;\,r;\,\nu)\in\Bbb N$ such that in
the $(m_1,\dots,m_N;\,r;\,\nu)$-billiard flow amongst any $N_0$ consecutive
collisions one can always find two neighboring ones separated in time by at
least $\epsilon_0$. Thus, for a phase point $x_{-t}$ at least
$\epsilon_0/2$-away from collisions, the norms $||\delta q_0||$ and
$\sqrt{||\delta q_0||^2+||\delta v_0||^2}$ are equivalent for all vectors
$(\delta q_0,\,\delta v_0)\in E^u(x_{-t})$, hence we immediately get

\medskip

\subheading{Corollary 3.14}
For any phase point $x_0\in\bold M\setminus\partial\bold M$ with a
non-singular backward trajectory $S^{(-\infty,0)}x_0$ and with infinitely many
consecutive, connected collision graphs on $S^{(-\infty,0)}x_0$, and for any
number $L>0$ one can always find a time $-t<0$ and a non-zero tangent vector
$(\delta q_0,\delta v_0)\in E^u(x_{-t})$ ($x_{-t}=S^{-t}x_0$) with

$$
\frac{||(\delta q_t,\,\delta v_t)||}{||(\delta q_0,\,\delta v_0)||}>L,
$$
where 
$(\delta q_t,\delta v_t)=DS^t(\delta q_0,\delta v_0)\in E^u(x_0)$. \qed

\medskip

The time-reversal dual of the previous result is immediately obtained by
replacing the phase point $x_0=(q_0,\,v_0)$ with $-x_0=(q_0,\,-v_0)$, the
backward orbit with the forward orbit, and the unstable vectors with the
stable ones. We formulate this dual as our

\medskip

\subheading{Theorem} 
For any phase point $x_0\in\bold M\setminus\partial\bold M$ with
a non-singular forward orbit $S^{(0,\infty)}x_0$ and with infinitely many
consecutive, connected collision graphs on $S^{(0,\infty)}x_0$, and for any
number $L>0$ one can always find a time $t>0$ and a non-zero tangent vector
$(\delta q_0,\delta v_0)\in E^s(x_{0})$ with

$$
\frac{||(\delta q_t,\,\delta v_t)||}{||(\delta q_0,\,\delta v_0)||}<L^{-1},
$$
where 
$(\delta q_t,\delta v_t)=DS^t(\delta q_0,\delta v_0)\in E^s(x_t)$,
$x_t=S^tx_0$. \qed


\bigskip \bigskip

\Refs

\widestnumber\key{B-Ch-Sz-T(2002)}

\ref\key B-Ch-Sz-T(2002)
\by P. B\'alint, N. Chernov, D. Sz\'asz, I. P. T\'oth
\paper Multi\-dimensional semi\-dispersing billiards: 
singularities and the fundamental theorem
\jour Ann. Henri Poincar\'e
\vol 3, No. 3 (2002)
\pages 451--482
\endref

\ref\key B-F-K(1998)
\by D. Burago, S. Ferleger, A. Kononenko
\paper Uniform estimates on the number of collisions in semi-dispersing
billiards
\jour Annals of Mathematics
\vol 147 (1998)
\pages 695-708
\endref

\ref\key B-R(1997)
\by L. A. Bunimovich, J. Rehacek
\paper Nowhere Dispersing 3D Billiards with Nonvanishing Lyapunov Exponents
\jour Commun. Math. Phys.
\vol 189 (1997), no. 3
\pages 729--757
\endref

\ref\key B-R(1998)
\by L. A. Bunimovich, J. Rehacek
\paper How High-Dimensional Stadia Look Like
\jour Commun. Math. Phys.
\vol 197 (1998), no. 2
\pages 277--301
\endref

\ref\key Ch(1982)
\by N. I. Chernov
\paper Construction of transverse fiberings in multidimensional 
semi-dispersed billiards
\jour Functional Anal. Appl.
\vol 16 (1982), no. 4
\pages 270--280
\endref

\ref\key Ch(1994)
\by N. I. Chernov
\paper Statistical Properties of the Periodic Lorentz Gas.
Multidimensional Case
\jour Journal of Statistical Physics
\vol 74, Nos. 1/2 (1994)
\pages 11-54
\endref

\ref\key Ch-Sim(2006)
\by N. I. Chernov, N. Sim\'anyi
\paper Flow-invariant hypersurfaces in semi-dispersing billiards
\jour To appear in Annales Henri Poincar\'e, arxiv:math.DS/0603360
\year 2006
\endref

\ref\key G(1981)
\by     G. Galperin
\paper On systems of locally interacting and repelling particles moving in
space
\jour Trudy MMO
\vol 43 (1981)
\pages 142-196
\endref

\ref\key K-B(1994)
\by A. Katok, K. Burns
\paper Infinitesimal Lyapunov functions, invariant cone families and
stochastic properties of smooth dynamical systems
\jour Ergodic Theory Dyn. Syst.
\vol 14, No. 4
\year 1994
\pages 757-785
\endref

\ref\key K-S-Sz(1990)-I
\by A. Kr\'amli, N. Sim\'anyi, D. Sz\'asz
\paper A ``Transversal'' Fundamental Theorem for Semi-Dis\-pers\-ing Billiards
\jour Commun. Math. Phys.
\vol 129 (1990)
\pages 535--560
\endref

\ref\key L-W(1995)
\by C. Liverani, M. Wojtkowski
\paper Ergodicity in Hamiltonian systems
\jour Dynamics Reported
\vol 4 (1995)
\pages 130-202, arXiv:math.DS/9210229
\endref

\ref\key Sim(1992)-I
\by N. Sim\'anyi
\paper The K-property of $N$ billiard balls I
\jour Invent. Math.
\vol 108 (1992)
\pages 521-548
\endref

\ref\key Sim(2002)
\by N. Sim\'anyi
\paper The Complete Hyperbolicity of Cylindric Billiards
\jour Ergodic Th. \& Dyn. Sys.
\vol 22 (2002)
\pages 281-302
\endref

\ref\key Sim(2003)
\by N. Sim\'anyi
\paper Proof of the Boltzmann-Sinai Ergodic Hypothesis for Typical Hard Disk 
Systems
\jour Inventiones Mathematicae
\vol 154, No. 1 (2003)
\pages 123-178
\endref

\ref\key Sim(2006)
\by N. Sim\'anyi
\paper Conditional Proof of the Boltzmann-Sinai Ergodic Hypothesis
(Assuming the Hyperbolicity of Typical Singular Orbits)
\jour Submitted for publication
\pages arXiv:math/0605358
\endref

\ref\key Sin(1979)
\by Ya. G. Sinai
\paper Development of Krylov's ideas. Afterword to N. S. Krylov's
``Works on the foundations of statistical physics''
\jour Princeton University Press, 1979
\endref

\ref\key S-Ch(1987)
\by Ya. G. Sinai, N.I. Chernov
\paper Ergodic properties of certain systems of 2--D discs and 3--D balls
\jour Russian Math. Surveys
\vol 42, No. 3 (1987)
\pages 181-207
\endref

\ref\key S-Sz(1999)
\by N. Sim\'anyi, D. Sz\'asz
\paper Hard ball systems are completely hyperbolic
\jour Annals of Mathematics
\vol 149 (1999)
\pages 35-96
\endref

\ref\key S-Sz(2000)
\by N. Sim\'anyi, D. Sz\'asz
\paper Non-integrability of Cylindric Billiards and
Transitive Lie Group Actions
\jour Ergod. Th. \& Dynam. Sys.
\vol 20 (2000)
\pages 593-610
\endref

\ref\key V(1979)
\by L. N. Vaserstein
\paper On Systems of Particles with Finite Range and/or Repulsive
Interactions
\jour Commun. Math. Phys.
\vol 69 (1979)
\pages 31-56
\endref

\endRefs
\bye